\documentclass[10pt, a4paper]{article}
\usepackage{mathtext}
\usepackage{mathrsfs}
\usepackage{amssymb}
\usepackage{amsfonts}
\usepackage{tensor}
\usepackage{amscd}
\usepackage{amsmath, amsthm}
\usepackage{mathtools}
\usepackage{hyperref}
\usepackage[nice]{nicefrac}
\usepackage[symbol]{footmisc}
\usepackage{empheq}

\begin{document}
\begin{center}\textbf{ON FAMILIES OF MONIC POLYNOMIALS}
\end{center}
\begin{center}
           \begin{footnotesize}
	DANIL KROTKOV
	\end{footnotesize}
\end{center}

\begin{footnotesize}
\noindent\textsc{Abstract}. In this paper we derive generalizations of different properties of monic polynomial families of binomial type, i.e. families of monic polynomials, for which the binomial theorem holds
$$
p_n(\alpha+\beta)=\sum_{k=0}^n \left(\vphantom{\bigg|}\genfrac{}{}{0pt}{0}{n}{k}\right) p_k(\alpha)p_{n-k}(\beta)
$$
Some trivial representations of general ''multiplication'' and ''derivative'' operators are derived. In addition we derive a formula for the logarithmic derivative of general monic polynomial $p_n(x)$ which reduces to the formula
$$
\frac{1}{n}\frac{p_n'(x)}{p_n(x)} =\left(x+\frac{1}{\varphi'(y)}\left(\frac{d}{dy}-n\mathrm{L}\right)\right)^{-1}\cdot\left.\frac{\varphi(y)}{y\varphi'(y)}~\right|_{y=0}
$$
derived by the author in binomial case, when the generating function of $p_n(x)$ equals to $e^{x\varphi(y)}$.
\end{footnotesize}\\

\noindent\underline{\textbf{Definition}}. We call a sequence of polynomials $\mathrm{P}=\{p_n(x)\}_{n\in\mathbb{N}_{0}}$ by a family of monic polynomials, if $\mathrm{deg}~p_n(x)=n$ and the leading coefficient of each polynomial is equal to $1$. Throughout this paper we assume that all polynomials and formal power series are defined over $\mathbb{C}$.\\

\noindent For given family of monic polynomials $\mathrm{P}$ there exists linear operator $\mathrm{G_{P}^{\vphantom{-1}}}$, which acts on polynomials by the rule
$$
\mathrm{G_{P}^{\vphantom{-1}}} \cdot x^n = p_n(x)
$$
\noindent This operator is invertible, since every monic polynomial may be expressed as a linear combination of $p_k(x)$ uniquely. In this paper we will be interested in representations of operators over polynomials as linear combinations of other operators. As an example, we write the simplest representation for the operator $\mathrm{G_{P}^{\vphantom{-1}}}$ and its inverse after the following definition.\\

\noindent \underline{\textbf{Definition}}. For given sequence $\{ a_n\}_{n\in\mathbb{N}_{0}}$ define an operator $a_\theta$: $a_\theta \cdot x^n = a_n x^n$. If $a_n=h(n)$ for some polynomial $h(x)$, then the following relation holds: $a_\theta=h(\theta)=h(x\mathrm{D})$, where $\mathrm{D}=\frac{d}{dx}$.\\

\noindent Suppose then that polynomial $p_n(x)$ has the following coefficients, and suppose that we have the following representation of monomial $x^n$ as linear combination of $p_k(x)$
$$
p_n(x)=\sum_{k=0}^{n}\binom{n}{k}~\xi_{k}^{n}x^{n-k}~~~~~~~~~~~~~~~~x^n=\sum_{k=0}^{n}\binom{n}{k}\mathbin{^*\!\xi}_{k}^{n}p_{n-k}(x)
$$
Then the following representations are valid (note that $\xi_{0}^{n}=\mathbin{^*\!\xi}_{0}^{n}=1$, since polynomials are monic)
$$
\mathrm{G_{P}^{\vphantom{-1}}}=\sum_{k=0}^\infty\frac{\mathrm{D}^k}{k!}\xi_{k}^{\theta}~~~~~~~~~~~~~~~~\mathrm{G}_{\mathrm{P}}^{-1}=\sum_{k=0}^\infty\frac{\mathrm{D}^k}{k!}\mathbin{^*\!\xi}_{k}^{\theta}
$$
It is natural then to consider the following ''multiplication'' and ''derivative'' operators
\begin{align*}
\mathrm{U_P}\coloneqq\mathrm{G_{P}^{\vphantom{-1}}}x\mathrm{G}_{\mathrm{P}}^{-1}~~~~~&\Rightarrow~~~~~\mathrm{U_P}\cdot p_n(x)=p_{n+1}(x)\\
\mathrm{D_P}\coloneqq\mathrm{G_{P}^{\vphantom{-1}}}\mathrm{D}\mathrm{G}_{\mathrm{P}}^{-1}~~~~~&\Rightarrow~~~~~\mathrm{D_P}\cdot p_n(x)=np_{n-1}(x)
\end{align*}
with an obvious relation
$$
\mathrm{D_P}\mathrm{U_P}-\mathrm{U_P}\mathrm{D_P}=1
$$
In addition, for given sequence $\{ a_n\}_{n\in\mathbb{N}_{0}}$ define an operator $a_{\mathrm{U_P}\mathrm{D_P}}^{\vphantom{-1}}\coloneqq \mathrm{G_{P}^{\vphantom{-1}}}a_{\theta}\mathrm{G}_{\mathrm{P}}^{-1}$. To move on we make the following observation. For given family $\mathrm{P}$ there exists a family of fomal power series $\mathrm{F}=\{f_n(y)\}_{n\in\mathbb{N}_{0}}$, s.t. $f_n(y) \in y^n+y^{n+1}\mathbb{C}[[y]]$, and such that the following relation holds
$$
e^{xy}=\sum_{n=0}^\infty \frac{p_n(x)f_n(y)}{n!}
$$
The series $f_n$ may be written down explicitly by the use of polynomials, which belong to the dual family $\mathrm{P}^*=\{p_n^*(x)\}_{n\in\mathbb{N}_{0}}$ of polynomials defined by the relation $p_n^*(x)\coloneqq\mathrm{G}_{\mathrm{P}}^{-1}\cdot x^n$.
\newpage\noindent
$$
x^n=\sum_{k=0}^{n}\binom{n}{k}\mathbin{^*\!\xi}_{k}^{n}p_{n-k}(x)~~~ \Rightarrow ~~~p_n^*(x)=\sum_{k=0}^{n}\binom{n}{k}\mathbin{^*\!\xi}_{k}^{n}x^{n-k}~~~ \Rightarrow~~~ f_n(y) = \sum_{k=0}^\infty \frac{y^{n+k}}{k!}\mathbin{^*\!\xi}_{k}^{n+k}
$$
Note that in case of binomial monic family there always exists such a series $f(y)\in y+y^2\mathbb{C}[[y]]$ that $f_n(y)=f(y)^n$. In this paper we are interested in different representations of the operators $\mathrm{U_P}$, $\mathrm{D_P}$ and in mixed relations between them and standard operators $x$, $\mathrm{D}$. For example, in binomial case $\mathrm{P}^{\mathrm{binom}}$ with $\mathrm{F}^{\mathrm{binom}}=\{f(y)^n\}_{n\in\mathbb{N}_{0}}$, we have the relations (see \cite{Dlt})
$$
\mathrm{U_P}=x\frac{1}{f'(\mathrm{D})}~~~~~~~~~~~~~~~~\mathrm{D_P}=f(\mathrm{D})
$$
There is also an interesting case of general orthogonal polynomials when there exist sequences $\{a_n\}_{n\in\mathbb{N}_{0}}$, $\{b_n\}_{n\in\mathbb{N}_{0}}$ and there is a mixed relation (see \cite{Gdsl})
$$
\mathrm{U_P}=x-a_{\mathrm{U_P}\mathrm{D_P}}^{\vphantom{-1}}-\mathrm{D_P}b_{\mathrm{U_P}\mathrm{D_P}}^{\vphantom{-1}} ~~~\iff~~~ \mathrm{U_{P^{\text{*}}}}=x+a_\theta+\mathrm{D}b_\theta
$$
but this case will not be overviewed in this paper. For now we prove the following lemma to derive some basic formulae for $\mathrm{U_P}$, $\mathrm{D_P}$ in general case.\\

\noindent\underline{\textbf{Lemma}}. ~~$f_0(y)=1 ~~\iff~~ \forall n>0,~~ p_n(0)=0$.\\

\noindent~~\textit{Proof}. Suppose $\forall n>0~ p_n(0)=0$. Then
$$
e^{xy}=\sum_{n=0}^\infty \frac{p_n(x)f_n(y)}{n!} ~~~~\Rightarrow~~~~ 1=\sum_{n=0}^\infty \frac{p_n(0)f_n(y)}{n!}=f_0(y)
$$
~~Suppose $f_0(y)=1$. We have $f_n(y) \in y^n+y^{n+1}\mathbb{C}[[y]]$, and thus
$$
e^{xy}=1+\sum_{n=1}^\infty \frac{p_n(x)f_n(y)}{n!} ~~~~\Rightarrow~~~~ \forall k\geqslant 1,~~ x^k=\sum_{n=1}^k \frac{p_n(x)}{n!}f_n^{(k)}(0).
$$
~~That means $p_1(x)=x$ and $\forall k \geqslant 2$ we have the relation
$$
p_k(x)=\frac{p_k(x)}{k!}f_k^{(k)}(0)=x^k-\sum_{n=1}^{k-1}\frac{p_n(x)}{n!}f_n^{(k)}(0)
$$
~~Now $p_1(0)=0$, and hence by induction we have also $~\forall n\geqslant 2, ~~p_n(0)=0$. \qed\\

\noindent It then follows that the sequence $\widetilde{\mathrm{P}}=\{\nicefrac{p_{n+1}(x)}{x}\}_{n\in\mathbb{N}_{0}}$ is a family of monic polynomials whenever $f_0(y)=1$. Moreover, for any monic polynomial family $\mathrm{P}=\{p_n(x)\}_{n\in\mathbb{N}_{0}}$ there exists another family $\widetilde{\mathrm{P}}=\{\widetilde{p}_n(x)\}_{n\in\mathbb{N}_{0}}=\{\nicefrac{f_0(\mathrm{D})\cdot p_{n+1}(x)}{x}\}_{n\in\mathbb{N}_{0}}$. The latter holds, since we have
$$
f_0(y)e^{xy}=f_0(\mathrm{D})\cdot e^{xy}=\sum_{n=0}^\infty \frac{[f_0(\mathrm{D})\cdot p_n(x)]f_n(y)}{n!}~~~ \Rightarrow~~~ e^{xy}=\sum_{n=0}^\infty \frac{[f_0(\mathrm{D})\cdot p_n(x)]}{n!}\frac{f_n(y)}{f_0(y)}
$$
Hence, it follows by the lemma above that polynomials $f_0(\mathrm{D})\cdot p_n(x)$ vanish at $x=0$ for all $n\geqslant 1$. Acting with $\frac{\partial}{\partial y}$ on both sides of the last equality, we obtain
$$
e^{xy}=\sum_{n=0}^\infty \frac{1}{n!}\frac{f_0(\mathrm{D})\cdot p_{n+1}(x)}{x}\frac{1}{n+1}\left(\frac{f_{n+1}(y)}{f_0(y)}\right)'=\sum_{n=0}^\infty \frac{\widetilde{p}_n(x)\widetilde{f}_n(y)}{n!}
$$
\noindent\underline{\textbf{Lemma}}. Suppose we have two monic families $\mathrm{P}=\{p_n(x)\}_{n\in\mathbb{N}_{0}}$, $\mathrm{Q}=\{q_n(x)\}_{n\in\mathbb{N}_{0}}$ and corresponding families $\mathrm{F}=\{f_n(y)\}_{n\in\mathbb{N}_{0}}$, $\mathrm{E}=\{e_n(y)\}_{n\in\mathbb{N}_{0}}$
$$
e^{xy}=\sum_{n=0}^\infty \frac{p_n(x)f_n(y)}{n!}=\sum_{n=0}^\infty \frac{q_n(x)e_n(y)}{n!}
$$
Then the following formula holds
$$
\mathrm{G_{P}^{\vphantom{-1}}}\mathrm{G_{Q}^{-1}}=(x)\frac{p_{\mathrm{U_Q}\mathrm{D_Q}}^{\vphantom{-1}}}{q_{\mathrm{U_Q}\mathrm{D_Q}}^{\vphantom{-1}}}=\frac{e_{\mathrm{U_P}\mathrm{D_P}}^{\vphantom{-1}}}{f_{\mathrm{U_P}\mathrm{D_P}}^{\vphantom{-1}}}(\mathrm{D})
$$
\newpage\noindent
~~\textit{Proof}. The first equality of operators is obvious
$$
(x)\frac{p_{\mathrm{U_Q}\mathrm{D_Q}}^{\vphantom{-1}}}{q_{\mathrm{U_Q}\mathrm{D_Q}}^{\vphantom{-1}}}\cdot q_n(x)=(x)\frac{p_n}{q_n}~q_n(x)=p_n(x)=\mathrm{G_{P}^{\vphantom{-1}}}\mathrm{G_{Q}^{-1}}\cdot q_n(x)
$$
~~hence these operators coincide on any polynomial. In case of the second equality we have
\begin{align*}
\frac{e_{\mathrm{U_P}\mathrm{D_P}}^{\vphantom{-1}}}{f_{\mathrm{U_P}\mathrm{D_P}}^{\vphantom{-1}}}(\mathrm{D})\cdot e^{xy}&=\frac{e_{\mathrm{U_P}\mathrm{D_P}}^{\vphantom{-1}}}{f_{\mathrm{U_P}\mathrm{D_P}}^{\vphantom{-1}}}(y)\cdot e^{xy}=\frac{e_{\mathrm{U_P}\mathrm{D_P}}^{\vphantom{-1}}}{f_{\mathrm{U_P}\mathrm{D_P}}^{\vphantom{-1}}}(y)\cdot\sum_{n=0}^\infty \frac{p_n(x)f_n(y)}{n!}=\sum_{n=0}^\infty \frac{e_n}{f_n}(y)\frac{p_n(x)f_n(y)}{n!}=\\
\tag*{\qed}&=\sum_{n=0}^\infty\frac{p_n(x)e_n(y)}{n!}=\mathrm{G_{P}^{\vphantom{-1}}}\mathrm{G_{Q}^{-1}}\cdot\sum_{n=0}^\infty\frac{q_n(x)e_n(y)}{n!}=\mathrm{G_{P}^{\vphantom{-1}}}\mathrm{G_{Q}^{-1}}\cdot e^{xy}
\end{align*}
Now we are ready to derive the following representations for $\mathrm{U_P}$, $\mathrm{D_P}$.\\

\noindent\underline{\textbf{Proposition}}. The operators $\mathrm{U_P}$, $\mathrm{D_P}$ have the following mixed representations
$$
\mathrm{D_P}=\frac{f_{\mathrm{U_P}\mathrm{D_P}+1}^{\vphantom{-1}}}{~f_{\mathrm{U_P}\mathrm{D_P}}^{\vphantom{-1}}}(\mathrm{D});~~~~~~~~~~\mathrm{U_P}=f_0(\mathrm{D})^{-1}x\mathrm{G_{\widetilde{P}}}\mathrm{G_{P}^{-1}} ~~~\Rightarrow~~~ x=f_0(\mathrm{D})\mathrm{U_P}\frac{\widetilde{f}_{\mathrm{U_P}\mathrm{D_P}}^{\vphantom{-1}}}{f_{\mathrm{U_P}\mathrm{D_P}}^{\vphantom{-1}}}(\mathrm{D})
$$
~~\textit{Proof}. We have
\begin{align*}
\mathrm{D_P}\cdot e^{xy}&=\sum_{n=1}^\infty \frac{np_{n-1}(x)f_n(y)}{n!}=\sum_{n=0}^\infty \frac{p_n(x)f_n(y)}{n!}\frac{f_{n+1}}{f_n}(y)=\sum_{n=0}^\infty \frac{f_{\mathrm{U_P}\mathrm{D_P}+1}^{\vphantom{-1}}}{~f_{\mathrm{U_P}\mathrm{D_P}}^{\vphantom{-1}}}(y)\cdot \frac{p_n(x)f_n(y)}{n!}=\\
&=\frac{f_{\mathrm{U_P}\mathrm{D_P}+1}^{\vphantom{-1}}}{~f_{\mathrm{U_P}\mathrm{D_P}}^{\vphantom{-1}}}(y)\cdot e^{xy}=\frac{f_{\mathrm{U_P}\mathrm{D_P}+1}^{\vphantom{-1}}}{~f_{\mathrm{U_P}\mathrm{D_P}}^{\vphantom{-1}}}(\mathrm{D})\cdot e^{xy}
\end{align*}
~~and
\begin{align*}
&\mathrm{U_P} \cdot p_n(x)=p_{n+1}(x)=f_0(\mathrm{D})^{-1}x\frac{f_0(\mathrm{D})\cdot p_{n+1}(x)}{x}=f_0(\mathrm{D})^{-1}x\mathrm{G_{\widetilde{P}}}\mathrm{G_{P}^{-1}}\cdot p_n(x) ~~\Rightarrow~~~~~\\
\tag*{\qed}&~~~~\mathrm{U_P}=f_0(\mathrm{D})^{-1}x\mathrm{G_{\widetilde{P}}}\mathrm{G_{P}^{-1}} ~~\Rightarrow~~ x=f_0(\mathrm{D})\mathrm{U_P}\mathrm{G_{P}^{\vphantom{p-1}}}\mathrm{G_{\widetilde{P}}^{\vphantom{p}-1}}=f_0(\mathrm{D})\mathrm{U_P}\frac{\widetilde{f}_{\mathrm{U_P}\mathrm{D_P}}^{\vphantom{-1}}}{f_{\mathrm{U_P}\mathrm{D_P}}^{\vphantom{-1}}}(\mathrm{D})
\end{align*}
Note that these trivial identities may actually be helpful, since they sometimes allow us to understand $\mathrm{U_P}$ or $\mathrm{D_P}$ without understanding $\mathrm{G_{P}^{\vphantom{-1}}}$ or its inverse. For example, in the simplest binomial case we have $f_{n+1}(y)f_{n}(y)^{-1}=f(y)$ for some $f(y)$, hence we obtain that indeed $\mathrm{D_P}=f(\mathrm{D})$. And also $\widetilde{f}_n=(1+n)^{-1}(f^{n+1})'=f'f^n$, thus $\widetilde{f}_{n}(y)f_{n}(y)^{-1}=f'(y)$ and hence $x=\mathrm{U_P}f'(\mathrm{D})$, as expected. In general such an approach gives mixed identities. For example, consider the case $f_n(y)=y^n \exp(n(n-1)y/2)$. Then the following holds
$$
\frac{f_{n+1}(y)}{f_n(y)}=ye^{ny} ~~~~\Rightarrow ~~~~\mathrm{D_P}=\sum_{k=0}^\infty \frac{(\mathrm{U_P}\mathrm{D_P})^k \mathrm{D}^{k+1}}{k!}
$$
Using the same approach one can derive similar identities, such as
$$
\mathrm{D}_{\mathrm{P}}^k=\frac{f_{\mathrm{U_P}\mathrm{D_P}+k}^{\vphantom{-1}}}{~f_{\mathrm{U_P}\mathrm{D_P}}^{\vphantom{-1}}}(\mathrm{D})
$$
Additional identities may be obtained by the use of conjugation with $\mathrm{G}(\cdot)\mathrm{G}^{-1}$, for example
$$
\mathrm{D_P}=\frac{f_{\mathrm{U_P}\mathrm{D_P}+1}^{\vphantom{-1}}}{~f_{\mathrm{U_P}\mathrm{D_P}}^{\vphantom{-1}}}(\mathrm{D})~~~~\Rightarrow ~~~~\mathrm{D}=\frac{f_{\theta+1}}{f_{\theta}}(\mathrm{D}_{\mathrm{P^{\text{*}}}})=\frac{f_{\theta+1}^*}{f_{\theta}^*}(\mathrm{D}_{\mathrm{P}})
$$
The latter proposition may be naturally generalized to an arbitrary operator on polynomials. In order to do that we have to introduce the following definition.
\newpage\noindent
\underline{\textbf{Definition}}. For an arbitrary operator $\mathrm{T}_x$, which acts on $\mathbb{C}[x]$, we define an operator $\mathrm{\overline{T}}_y$, which acts on $\mathbb{C}[[y]]$, such that
$$
\mathrm{T}_x\cdot e^{xy}=\mathrm{\overline{T}}_y \cdot e^{xy}
$$
Such an operator always exists and is fully determined by its images on monomials, since for an arbitrary $\ell(y) \in \mathbb{C}[[y]]$ holds
$$
\mathrm{\overline{T}}_y \cdot \ell(y)= \mathrm{\overline{T}}_y \cdot \ell(y)e^{xy}\big|_{x=0}=\mathrm{\overline{T}}_y \ell(\mathrm{D})\cdot e^{xy}\big|_{x=0}=\ell(\mathrm{D})\mathrm{\overline{T}}_y\cdot e^{xy}\big|_{x=0}=\sum_{k=0}^\infty \frac{y^k}{k!} [\ell(\mathrm{D})\mathrm{T}_x \cdot x^n]\big|_{x=0}
$$
Since $\mathrm{T}_x$ sends polynomials to polynomials, the action of $\ell(\mathrm{D})$ is well-defined and thus the action of $\mathrm{\overline{T}}_y$ is well-defined too. Note that the action of $\mathrm{T}_x$ can be reconstructed from the action of $\mathrm{\overline{T}}_y$. We denote $\mathrm{\overline{U}_P}$ by $\mathfrak{d}_{\mathrm{P}}$ and $\mathrm{\overline{D}_P}$ by $\mathfrak{u}_{\mathrm{P}}$. Note that in case of binomial monic family $\mathrm{P^{binom}}$ with $\mathrm{F^{binom}}=\{f(y)^n\}_{n \in \mathbb{N}_0}$ we have $\mathrm{D_P}=f(\mathrm{D})$, hence $\mathfrak{u}_{\mathrm{P}}=f(y)$, and $\mathrm{U_P}=xf'(\mathrm{D})^{-1}$, hence $\mathfrak{d}_{\mathrm{P}}=\frac{d}{df}$. Now we are ready to formulate the following theorem on representations of operators over $\mathbb{C}[x]$. For an arbitrary monic polynomial family $\mathrm{P}$ consider the generating function
$$
Q_y^{\mathrm{P}}(x)\coloneqq \sum_{n=0}^\infty \frac{p_n(x)y^n}{n!}=\sum_{n=0}^\infty \frac{x^n f_n^*(y)}{n!} ~~~~\Rightarrow~~~~\mathrm{D_P}\cdot Q_y^{\mathrm{P}}(x)=yQ_y^{\mathrm{P}}(x); ~~~~~ \mathfrak{d}_{\mathrm{P}^{\text{*}}}\cdot Q_y^{\mathrm{P}}(x)=xQ_y^{\mathrm{P}}(x)
$$
\noindent\underline{\textbf{Theorem}}. An arbitrary operator $\mathrm{T}_x$ on $\mathbb{C}[x]$ has the following representations
\begin{align*}
\mathrm{T}_x&=(x)\frac{\mathrm{T}\cdot Q_y^{\mathrm{P}}}{Q_y^{\mathrm{P}}}\bigg|_{y=\mathrm{D_P}}=(x)\frac{\mathrm{T}\cdot p_n}{p_n}\bigg|_{n=\mathrm{U_P D_P}}=(\mathrm{U_P})\frac{\mathrm{\overline{T}}\cdot Q_y^{\mathrm{P}^{*}}}{Q_y^{\mathrm{P}^{\text{*}}}}\bigg|_{y=\mathrm{D}}=\frac{\mathrm{\overline{T}}\cdot f_n}{f_n}\bigg|_{n=\mathrm{U_P D_P}}(\mathrm{D})
\end{align*}
~~\textit{Proof}. The following holds for an arbitrary $t$, hence it holds for an arbitrary polynomial in $x$
$$
(x)\frac{\mathrm{T}\cdot Q_y^{\mathrm{P}}}{Q_y^{\mathrm{P}}}\bigg|_{y=\mathrm{D_P}}\cdot Q_t^{\mathrm{P}}(x)=(x)\frac{\mathrm{T}\cdot Q_y^{\mathrm{P}}}{Q_y^{\mathrm{P}}}\bigg|_{y=t} \cdot Q_t^{\mathrm{P}}(x)=\frac{\mathrm{T}\cdot Q_t^{\mathrm{P}}(x)}{Q_t^{\mathrm{P}}(x)} Q_t^{\mathrm{P}}(x)=\mathrm{T}\cdot Q_t^{\mathrm{P}}(x)
$$
~~The second equality is trivial, since
$$
(x)\frac{\mathrm{T}\cdot p_n}{p_n}\bigg|_{n=\mathrm{U_P D_P}}\cdot p_m(x)=(x)\frac{\mathrm{T}\cdot p_n}{p_n}\bigg|_{n=m}\cdot p_m(x)=\frac{\mathrm{T}\cdot p_m(x)}{p_m(x)}p_m(x)=\mathrm{T}\cdot p_m(x)
$$
~~In the third case we have
$$
\overline{\left((\mathrm{U_P})\frac{\mathrm{\overline{T}}\cdot Q_y^{\mathrm{P}^{*}}}{Q_y^{\mathrm{P}^{\text{*}}}}\bigg|_{y=\mathrm{D}}\right)}=\frac{\mathrm{\overline{T}}\cdot Q_y^{\mathrm{P}^{*}}}{Q_y^{\mathrm{P}^{\text{*}}}}(\mathfrak{d}_{\mathrm{P}})
$$
~~Now since $\mathfrak{d}_{\mathrm{P}^{\text{*}}}\cdot Q_y^{\mathrm{P}}(t)=tQ_y^{\mathrm{P}}(t)$, the following holds
$$
\frac{\mathrm{\overline{T}}\cdot Q_y^{\mathrm{P}^{*}}}{Q_y^{\mathrm{P}^{\text{*}}}}(\mathfrak{d}_{\mathrm{P}})\cdot Q_y^{\mathrm{P}^{*}}(t)=\frac{\mathrm{\overline{T}}\cdot Q_y^{\mathrm{P}^{*}}}{Q_y^{\mathrm{P}^{\text{*}}}}(t)\cdot Q_y^{\mathrm{P}^{*}}(t)=\mathrm{\overline{T}}\cdot Q_y^{\mathrm{P}^{*}}(t)
$$
~~Hence
$$
\mathrm{\overline{T}}=\overline{\left((\mathrm{U_P})\frac{\mathrm{\overline{T}}\cdot Q_y^{\mathrm{P}^{*}}}{Q_y^{\mathrm{P}^{\text{*}}}}\bigg|_{y=\mathrm{D}}\right)} ~~~~~~~~\Rightarrow~~~~~~~~ \mathrm{T}=(\mathrm{U_P})\frac{\mathrm{\overline{T}}\cdot Q_y^{\mathrm{P}^{*}}}{Q_y^{\mathrm{P}^{\text{*}}}}\bigg|_{y=\mathrm{D}}
$$
~~Finally
\begin{align*}
\frac{\mathrm{\overline{T}}\cdot f_n}{f_n}\bigg|_{n=\mathrm{U_P D_P}}(\mathrm{D})\cdot e^{xy} &= \frac{\mathrm{\overline{T}}\cdot f_n}{f_n}\bigg|_{n=\mathrm{U_P D_P}}(y)\cdot e^{xy}=\frac{\mathrm{\overline{T}}\cdot f_n}{f_n}\bigg|_{n=\mathrm{U_P D_P}}(y)\cdot \sum_{n=0}^\infty \frac{p_n(x)f_n(y)}{n!}=\\
&=\sum_{n=0}^\infty \frac{\mathrm{\overline{T}}\cdot f_n}{f_n}(y)\frac{p_n(x)f_n(y)}{n!}=\sum_{n=0}^\infty \frac{p_n(x)\mathrm{\overline{T}}\cdot f_n(y)}{n!}=\mathrm{\overline{T}}\cdot e^{xy}=\mathrm{T}\cdot e^{xy}
\end{align*}
\hfill\ensuremath$\blacksquare$
\newpage\noindent One may notice a sneaky moment in the latter reasoning. Of course, it works perfectly well, when one deals with operators on polynomials. But it is often natural to try to extend action of these operators to arbitrary series of the form $x^s \mathbb{C}[[x^{-1}]]$ for any $s\in \mathbb{C}$. And here one has to be careful with transformations of such a kind. Consider the simplest operator $a_{\theta}$. Now take the trivial monic family $\{x^n\}_{n\in\mathbb{N}_0}$ with $Q_y(x)=e^{xy}$. Then the latter theorem implies that
$$
a_{\theta}=(x)\frac{a_{\theta}\cdot e^{xy}}{e^{xy}}\bigg|_{y=\mathrm{D}}=\sum_{n=0}^\infty \frac{x^n\mathrm{D}^n}{n!}\sum_{k=0}^n \binom{n}{k}(-1)^{n-k}a_k=\sum_{n=0}^\infty \binom{\theta}{n}\sum_{k=0}^n \binom{n}{k}(-1)^{n-k}a_k
$$
Of course, the latter equality is true when one considers both sides as operators on $\mathbb{C}[x]$, since for any monomial the sum on the right side converges and the images coincide. But in case of extension to $x^s \mathbb{C}[[x^{-1}]]$, the operator on the left side can be extended arbitrarily, as soon as the generalized sequence $\{a_s\}_{s\in\mathbb{C}}$ is determined. On the other hand such an extension does not have to be determined by its subsequence $\{a_n\}_{n\in\mathbb{N}_0}$. Moreover, the series
$$
a_s \neq \sum_{n=0}^\infty \binom{s}{n}\sum_{k=0}^n \binom{n}{k}(-1)^{n-k}a_k
$$
does not have to converge in general. That means that theorem is true in case of operators on $x^s \mathbb{C}[[x^{-1}]]$ only for a specific class of operators (for example, it holds for some polynomial $a(s)=a_s$). Consider now the following corollary of the theorem.\\

\noindent\underline{\textbf{Corollary}}. For an arbitrary monic family $\mathrm{P}$ with generating function $Q_y^{\mathrm{P}}(x)=Q_y(x)$ and the dual generating function $Q_y^{\mathrm{P}^{*}}(x)=Q_y^*(x)$, the operator $\theta$ has the following representations
$$
\theta=x(x)\frac{Q'_{\mathrm{D_P}}}{Q_{\mathrm{D_P}}}=(\mathrm{U_P})\left(y\frac{\partial}{\partial y}\ln Q_y^*\right)\bigg|_{y=\mathrm{D}}=x(x)\frac{p'_{\mathrm{U_P D_P}}}{p_{\mathrm{U_P D_P}}}=\frac{f'_{\mathrm{U_P D_P}}}{f_{\mathrm{U_P D_P}}}(\mathrm{D})\mathrm{D}
$$
If a family $\mathrm{P}$ is binomial, then $Q_y^*(x)=\exp(xf(y))$, $Q_y(x)=\exp(x\varphi(y))$, where $\varphi(y)=f^{inv}(y)$ is a functional inverse to $f(y)$. In this case the statement above degenerates to the following equalities:
\begin{align*}
\frac{Q'_y}{Q_y}(x)=\varphi(y) ~~~~~\Rightarrow~~~~~ \theta=x&\varphi(\mathrm{D_P}),~\text{and that is indeed true, since}~ \mathrm{D_P}=f(\mathrm{D}) ~\Rightarrow~ \varphi(\mathrm{D_P})=\mathrm{D}\\
\left(y\frac{\partial}{\partial y}\ln Q_y^*\right)(x)=xyf'(y) ~~~~~&\Rightarrow~~~~~\theta=\mathrm{U_P D}f'(\mathrm{D}),~\text{and that is indeed true, since}~\mathrm{U_P}=x\frac{1}{f'(\mathrm{D})}\\
f_n(y)=f^n(y), ~~&\text{hence}~~ \frac{f'_n}{f_n}(y)=n\frac{f'(y)}{f(y)} ~~~~~\Rightarrow~~~~~ \theta=\mathrm{U_P D_P} \frac{f'(\mathrm{D})\mathrm{D}}{f(\mathrm{D})}
\end{align*}
In other cases the corollary implies different mixed identities. Going back to our example $f_n(y)=y^n\exp(n(n-1)y/2)$, we obtain in this case
$$
\frac{f'_n}{f_n}(y)y=n+\frac{n(n-1)}{2}y ~~~~~\Rightarrow~~~~~ \theta=\mathrm{U_P D_P}+\frac{1}{2}\mathrm{U_P^2}\mathrm{D_P^2}\mathrm{D}
$$
These identities may actually be helpful, if one wants to understand the behaviour of monic polynomials. The equality of operators
$$
x(x)\frac{Q'_{\mathrm{D_P}}}{Q_{\mathrm{D_P}}}=x(x)\frac{p'_{\mathrm{U_P D_P}}}{p_{\mathrm{U_P D_P}}}
$$
actually tells us how to calculate the expansion of the logarithmic derivative in a natural way, when the generating function is known. The main problem here is to understand, how to construct expansion of an operator of the form
$$
\sum_{k=0}^\infty x^{-k}\mu_k^{\mathrm{U_P D_P}}
$$
from the expansion of the form
$$
\sum_{k=0}^\infty \ell_k(x)\mathrm{D}_{\mathrm{P}}^k
$$
\newpage\noindent Note that one can easily do that formally using symbolic calculation. We denote $\frac{d}{dy}=\mathfrak{d}$
\begin{align*}
\sum_{k=0}^\infty \ell_k(x)\mathrm{D}_{\mathrm{P}}^k&=\ell_0(x)+\sum_{k=1}^\infty \ell_k(x)\mathrm{D}_{\mathrm{P}}^k=\ell_0(x)+\sum_{k=1}^\infty \ell_k(x)\mathrm{D}_{\mathrm{P}}^{k-1}\mathrm{U}_{\mathrm{P}}^{-1}(\mathrm{U_P D_P})=\\
&=\ell_0(x)+\sum_{k=1}^\infty \ell_k(x)(x)\frac{\mathrm{D}_{\mathrm{P}}^{k-1}\mathrm{U}_{\mathrm{P}}^{-1} \cdot Q_y}{Q_y}\bigg|_{y=\mathrm{D_P}}(\mathrm{U_P D_P})=\\
&=\ell_0(x)+\sum_{k=1}^\infty \ell_k(x)(x)\frac{\mathfrak{d}^{-1} y^{k-1} \cdot Q_y}{Q_y}\bigg|_{y=\mathrm{D_P}}(\mathrm{U_P D_P})
\end{align*}
Now introduce the 0-derivative $\mathrm{L}: \mathrm{L}\cdot g(y)=\frac{g(y)-g(0)}{y}=(1+\theta)^{-1}\mathfrak{d}\cdot g(y)$. Then if one defines two auxiliary operators on $\mathbb{C}[[x,y]]$
$$
\mathrm{A}F(x,y)=F(x,0);~~~~~~~~~~~\mathrm{B}F(x,y)=Q_y(x)^{-1}\mathfrak{d}^{-1}Q_y(x)\mathrm{L}F(x,y)
$$
one formally obtains
$$
F(x,\mathrm{D_P})=(\mathrm{A}F)(x,\mathrm{D_P})+(\mathrm{B}F)(x,\mathrm{D_P})\mathrm{U_P D_P}
$$
Continuing this procedure, we obtain that the following expansion is valid
$$
F(x,\mathrm{D_P})=\sum_{k=0}^\infty (\mathrm{B}^k\cdot F(x,y))\bigg|_{y=0} (\mathrm{U_P D_P})^k
$$
Now act with this operator on $p_n(x)$ and notice the geometric progression to formally obtain
\begin{align*}
\frac{F(x,\mathrm{D_P})\cdot p_n(x)}{p_n(x)}&=(1-nQ_y(x)^{-1}\mathfrak{d}^{-1}Q_y(x)\mathrm{L})^{-1}\cdot F(x,y)\bigg|_{y=0}=\\
&=(Q_y(x)^{-1}\mathfrak{d}Q_y(x)-n\mathrm{L})^{-1}Q_y(x)^{-1}\mathfrak{d}Q_y(x)\cdot F(x,y)\bigg|_{y=0}
\end{align*}
The main problem in the latter reasoning is that we somehow defined the action of $\mathfrak{d}^{-1}$. Another problem is that we speculated that one can actually expand the fraction $(F(x,\mathrm{D_P})\cdot p_n(x))/p_n(x)$ into geometric progression as a function of $n$. So the latter formula may actually be true, but to prove it we have to be more accurate and more specific. One may actually notice the third problem here: how do we have to interpret the action of $(Q_y(x)^{-1}\mathfrak{d}Q_y(x)-n\mathrm{L})^{-1}$?\\

\noindent\underline{\textbf{Theorem}}. Suppose we have a monic family $\mathrm{P}$ such that $p_n(0)=0$ for all $n\geqslant 1$. Consider the generating function $Q_y(x)$ and define $Q'_y(x)\coloneqq\partial_x Q_y(x)=\mathrm{D}\cdot Q_y(x)$, $\dot Q_y(x)\coloneqq\partial_y Q_y(x)=\mathfrak{d}\cdot Q_y(x)$. Then the following formula holds
$$
\boxed{~~\frac{x}{n}\frac{p'_n(x)}{p_n(x)}=\left(1+\frac{Q_y(x)}{\dot Q_y(x)}(\partial_y-n\mathrm{L})\right)^{-1} \cdot \frac{x}{y}\frac{Q'_y(x)}{\dot Q_y(x)}\biggr|_{y=0}~~}
$$
where the action of an operator $(1+W)^{-1}$ is interpreted as an action of $\hat 1-W+W^2-...$.\\

\noindent Before we provide the proof of the theorem, we note first that the symbolic calculation above is not as useless, as it seems to be, since there is a specific case, when it works perfectly well, and that is the binomial one. This case allows us to avoid the flow with interpretation of the action of $\mathfrak{d}^{-1}$. It is still not well-defined, but the action of $Q_y(x)^{-1}\mathfrak{d}^{-1}Q_y(x)$ is well-defined in this case: since $Q_y(x)=\exp(x\varphi(y))$, we have $Q_y(x)^{-1}\mathfrak{d}^{-1}Q_y(x)=x^{-1}(1+x^{-1}\frac{d}{d\varphi})^{-1}\varphi'(y)^{-1}$. The second problem is not actually a problem too, since in binomial case we have a natural generalization of $p_n(x)$ to an arbitrary complex index, because the coefficients in the expansion $p_s(x)\in x^{s}+x^{s-1}\mathbb{C}[[x^{-1}]]$ are polynomials in $s$. We are now going to use the following fact about polynomials of binomial type and their continuations, derived by the author in the previous paper.
\newpage\noindent\underline{\textbf{Proposition}}. For an arbitrary operator $\mathrm{T}_{\alpha}$ and binomial family $\mathrm{P^{binom}}=\{p_n(\alpha)\}_{n\in\mathbb{N}_0}$ the following formula holds
$$
\sum_{n=0}^\infty \frac{p_n(\alpha)y^n}{n!}=e^{\alpha \varphi(y)} ~~\Rightarrow~~\frac{\mathrm{T}\cdot p_{n-1}(\alpha)}{p_n(\alpha)} =\left(\alpha+\frac{1}{\varphi'(y)}\left(\partial_y-n\mathrm{L}\right)\right)^{-1}\frac{1}{\varphi'(y)}e^{-\alpha\varphi(y)}\mathrm{T}\cdot e^{\alpha \varphi(y)}\biggr|_{y=0}
$$
~~\textit{Proof}. See formula (1.8) in \cite{LB}. \qed\\

\noindent We are now ready to prove the theorem. For monic family $\mathrm{P}$ such that $p_n(0)=0$ for all $n\geqslant 1$ consider the generating function $Q_y(x)$ and define
$$
\Phi_x(y)\coloneqq \frac{1}{x}\ln Q_y(x) \in y+y^2\mathbb{C}[x][[y]]
$$
Define also $f_x(y): f_x(\Phi_x(y))=y$. Consider now binomial family $\{\nu_n^x(\alpha)\}_{n\in\mathbb{N}_0}$ with generating function
$$
e^{\alpha \Phi_x(y)} = \sum_{n=0}^\infty \frac{\nu_n^x(\alpha) y^n}{n!}=(Q_y(x))^{\tfrac{\alpha}{x}}
$$
Consider also the following family of polynomials
$$
Q'_y(x)Q_y(x)^{-1}e^{\alpha \Phi_x(y)}=\sum_{n=0}^\infty \frac{\tau_n^x(\alpha) y^n}{n!}=[\partial_x x\cdot \Phi_x(y)]e^{\alpha \Phi_x(y)}
$$
Now notice that $\Phi_x(y) \in y+y^2\mathbb{C}[x][[y]] \Rightarrow [\partial_x x\cdot \Phi_x(y)] \in y+y^2\mathbb{C}[x][[y]]$, hence we have that $[\partial_x x\cdot \Phi_x(y)]/y \in 1+y\mathbb{C}[x][[y]]$. So we may define an operator $\mathrm{T}$ on polynomials in $\alpha$
$$
\mathrm{T}\coloneqq \frac{\partial_x x\cdot \Phi_x(t)}{t}\bigg|_{t=f_x(\partial/\partial\alpha)} ~~~~~\Rightarrow~~~~~\mathrm{T}\cdot e^{\alpha \Phi_x(y)}=\frac{1}{y}\frac{Q'_y(x)}{Q_y(x)} e^{\alpha \Phi_x(y)}
$$
By definition of $\nu_n^x(\alpha)$, $\tau_n^x(\alpha)$ and proposition, the following holds
\begin{align*}
\frac{1}{n}\frac{\tau_n^x(\alpha)}{\nu_n^x(\alpha)}&=\frac{1}{n}\frac{\mathrm{T}f_x(\partial/\partial\alpha)\cdot\nu_n^x(\alpha)}{\nu_n^x(\alpha)}=\frac{\mathrm{T}\cdot \nu_{n-1}^x(\alpha)}{\nu_n^x(\alpha)}=\\
&=\left(\alpha+\frac{x Q_y(x)}{\dot Q_y(x)}\left(\partial_y-n\mathrm{L}\right)\right)^{-1}\frac{x Q_y(x)}{\dot Q_y(x)}e^{-\alpha \Phi_x(y)}\mathrm{T}\cdot e^{\alpha \Phi_x(y)}\biggr|_{y=0}=\\
&=\left(\alpha+\frac{x Q_y(x)}{\dot Q_y(x)}\left(\partial_y-n\mathrm{L}\right)\right)^{-1}\cdot \frac{x Q_y(x)}{\dot Q_y(x)}\frac{1}{y}\frac{Q'_y(x)}{Q_y(x)}\biggr|_{y=0}=\\
&=\left(\alpha+\frac{x Q_y(x)}{\dot Q_y(x)}\left(\partial_y-n\mathrm{L}\right)\right)^{-1}\cdot \frac{x}{y}\frac{Q'_y(x)}{\dot Q_y(x)}\biggr|_{y=0}
\end{align*}
The latter equality essentially means that we have the following expansion
$$
\frac{\alpha}{n}\frac{\tau_n^x(\alpha)}{\nu_n^x(\alpha)}=\sum_{k=0}^\infty (-\alpha)^{-k} \left(\frac{x Q_y(x)}{\dot Q_y(x)}\left(\partial_y-n\mathrm{L}\right)\right)^k \cdot \frac{x}{y}\frac{Q'_y(x)}{\dot Q_y(x)}\biggr|_{y=0} = \sum_{k=0}^\infty  (-\alpha)^{-k}\xi_k(x)
$$
Hence, to prove the theorem it is enough to check that we can set $\alpha=x$, since $\tau_n^x(x)=p'_n(x)$ and $\nu_n^x(x)=p_n(x)$. For that purpose note that the expansion of $Q_y(x)e^{-xy}$ has polynomials of degree $k-1$ in $x$ as coefficients of $y^k$, since polynomials $p_n(x)$ are monic. The latter means that the expansion of $\ln Q_y(x)-xy$ has the same property. Also $p_{n>0}(0)=0$. Thus the following holds
\begin{align*}
\frac{1}{x}\frac{\dot Q_y(x)}{Q_y(x)}=1+\sum_{k=1}^\infty y^k h_k(x) ~~\text{\textit{\&}}~ \deg{h_k}=k-1;~~~~\frac{1}{y}\frac{Q'_y(x)}{Q_y(x)}=1+\sum_{k=1}^\infty y^k d_k(x) ~~\text{\textit{\&}}~ &\deg{d_k}=k-1\\
\frac{x Q_y(x)}{\dot Q_y(x)}=1+\sum_{k=1}^\infty y^k \ell_k(x) ~~\text{\textit{\&}}~ \deg{\ell_k}=k-1;~~~~\frac{x}{y}\frac{Q'_y(x)}{\dot Q_y(x)}=1+\sum_{k=1}^\infty y^k m_k(x) ~~\text{\textit{\&}}~ &\deg{m_k}=k-1
\end{align*}
\newpage\noindent Now set in addition $\ell_0(x)=m_0(x)=1$. We have the following formula for $\xi_M(x)$
\begin{align*}
\xi_M(x)&=\left(\sum_{k=0}^\infty \ell_k(x)y^k(\partial_y-n\mathrm{L})\right)^M \cdot~ \sum_{k=0}^\infty m_k(x)y^k ~\biggr|_{y=0}=\\
&=\left(\sum_{k=0}^\infty \ell_k(x)y^k(\partial_y-n\mathrm{L})\right)^{M-1} \cdot~ \sum_{k=0}^\infty \ell_k(x)y^k\sum_{k=0}^\infty (k+1-n)m_{k+1}(x)y^k ~\biggr|_{y=0}=\\
&=\left(\sum_{k=0}^\infty \ell_k(x)y^k(\partial_y-n\mathrm{L})\right)^{M-1} \cdot ~\sum_{k=0}^\infty y^k \sum_{q_1=0}^{k} (q_1+1-n)m_{q_1+1}(x)\ell_{k-q_1}(x) ~\biggr|_{y=0}=...=\\
&=\sum_{k=0}^\infty y^k \sum_{q_M=0}^k\sum_{q_{M-1}=0}^{q_M+1}...\sum_{q_1=0}^{q_2+1}(q_M+1-n)\ell_{k-q_M}\left[\prod_{i=1}^{M-1}(q_i+1-n)\ell_{q_{i+1}-q_{i}+1}\right]m_{q_1+1}~\biggr|_{y=0}=\\
&=\sum_{q_M=0}^{0}\sum_{q_{M-1}=0}^{q_M+1}\sum_{q_{M-2}=0}^{q_{M-1}+1}...\sum_{q_1=0}^{q_2+1}  m_{q_1+1}(x)\prod_{i=1}^{M}(q_i+1-n) \prod_{i=1}^{M-1}\ell_{q_{i+1}-q_{i}+1}(x)
\end{align*}
Notice that in case $n=1$ all nontrivial $\xi_M(x)$ vanish, since they are proportional to $(1-n)$. Hence we may assume that $n >1$. Now take an arbitrary summand with fixed indices $(q_M, q_{M-1},...,q_1)$. For convenience write $\sigma$, if $q_{i-1}=q_i+1$ and $\varepsilon_k$, if $q_{i-1}=q_i-k$, $(k\geqslant 0)$. We thus obtain a word of the form
$$
\underbrace{\sigma\sigma...\sigma}_{A_1}\varepsilon_{i_1}\underbrace{\sigma\sigma...\sigma}_{A_2}\varepsilon_{i_2}...\underbrace{\sigma\sigma...\sigma}_{A_{\ell-1}}\varepsilon_{i_{\ell-1}}\underbrace{\sigma\sigma...\sigma}_{A_{\ell}}
$$
where $A_1+1+A_2+1+...+A_{\ell-1}+1+A_{\ell}=M-1$ and $i_1 \leqslant A_1$,~...,~$i_{\ell-1}\leqslant A_1-i_1+...+A_{\ell-1}$. This word represents the summand with indices $(q_M, q_{M-1},...,q_1)$ with the following $q_i$
\begin{align*}
&q_M=0; ~~~~~~~~~~~~~~~~~~~~~~~~~~~~~~~~~~q_{M-1}=1; ~~~~~~~~~~~~~~~~~~~~~...~~~~~~q_{M-A_1}=A_1;\\
&q_{M-A_1-1}=A_1-i_1; ~~~~~~~~~~~~~~~~~~q_{M-A_1-2}=A_1-i_1+1;~~~...~~~~~~q_{M-A_1-1-A_2}=A_1-i_1+A_2;\\
&q_{M-A_1-A_2-2}=A_1-i_1+A_2-i_2;~~~~~~~~~~~~~~~~~~~...~~~~~~~~~~~~~~~~~~~~~~~~~~~~~...\\
&q_{M-A_1-A_2-...-A_{\ell-1}-\ell+1}=A_1-i_1+...+A_{\ell-1}-i_{\ell-1};~~...~~~q_{1}=A_1-i_1+...+A_{\ell-1}-i_{\ell-1}+A_{\ell}.
\end{align*}
Hence the coefficient is equal to
\begin{align*}
\prod_{i=1}^{M}(q_i+1-n)&=(-1)^M\prod_{i=0}^{M-1}(n-1-q_{M-i})=(-1)^M\prod_{k=0}^{\ell-1}~\prod_{j=A_1+...+A_k+k}^{A_1+...+A_{k+1}+k}~(n-1-q_{M-j})=\\
&=(-1)^M\prod_{k=0}^{\ell-1}~\prod_{j=0}^{A_{k+1}}(n-1-j-(A_1-i_1)-...-(A_k-i_k))
\end{align*}
and the polynomial is
$$
(\ell_{i_1+1}\ell_{i_2+1}...\ell_{i_{\ell-1}+1})(x)~m_{A_1-i_1+...+A_{\ell-1}-i_{\ell-1}+A_{\ell}+1}(x)
$$
Now since $\deg \ell_{k+1}=k$ and $\deg m_{k+1}=k$, the degree of the summand with nonzero coefficient is in general equal to $\deg S_{(q_M,...,q_1)}=i_1+i_2+...+i_{\ell-1}+A_1-i_1+...+A_{\ell-1}-i_{\ell-1}+A_{\ell}=M-\ell$. We now write the following conditions for coefficient to be nonzero. For any $0\leqslant k \leqslant \ell-1$:
\begin{align*}
A_{k+1}+1\leqslant n-1 -\sum_{m=1}^{k}(A_m-i_m) ~~~~\Rightarrow~~~~A_{k+1}+1\leqslant n-1
\end{align*}
Summing up the resulting inequalities we obtain that $M \leqslant (n-1)\ell$. Hence the degree of any nonvanishing summand is bounded by
$$
\deg S_{(q_M,...,q_1)}=M-\ell\leqslant M\frac{n-2}{n-1} ~~~~~~\Rightarrow~~~~~ \deg \xi_M(x) \leqslant M\frac{n-2}{n-1}
$$
Thus the series is well-defined at $\alpha=x$ and the theorem is proved. \hfill\ensuremath\blacksquare
\newpage
\begin{center}
\textbf{Conclusion}
\end{center}
One can see from the theorem, that there is a way to find logarithmic derivative of $f_n(x)$ from the generating function, since every theorem has a dual one in this structure. Note only that we still have to set $f_0^{*}(y)=1$. The formula then is as follows
\begin{align*}
Q_y(x)=\sum_{n=0}^\infty \frac{p_n(x)y^n}{n!}=\sum_{n=0}^\infty \frac{x^n f_n^*(y)}{n!}~~~~&\Rightarrow\\
\Rightarrow~~~~~~\frac{y}{n}\frac{(f_n^*(y))'}{f_n^*(y)}&=\left(1+\frac{Q_y(x)}{Q'_y(x)}(\partial_x-n\mathrm{L}_x)\right)^{-1} \cdot \frac{y}{x}\frac{\dot Q_y(x)}{Q'_y(x)}\biggr|_{x=0}
\end{align*}
and the reasoning above still works to prove the convergence of this series. But for now we should check the lowest coefficient of formal power series instead of the highest coefficient of a polynomial. Thus, redefining ''$\deg$'' for $\mathbb{C}[[y]]$, one may obtain the following expansion
$$
\frac{y}{n}\frac{(f_n^*(y))'}{f_n^*(y)}=\sum_{k\geqslant 0}(-y)^{-k}\psi_k(y)
$$
where now $\deg \psi_M(y)=M+\ell\geqslant M\tfrac{n}{n-1}$, hence the series does not include any negative power of $y$ and any coefficient of $y^h$ is well-defined. In other words, it doesn't really matter, if the coefficients $g_n(x)$ in the expansion
$$
F(x,y)=\sum_{n=0}^\infty \frac{g_n(x)y^n}{n!}
$$
are monic polynomials, or the series of the form $x^n+x^{n+1}\mathbb{C}[[x]]$. The theorem provides correct formula for the logarithmic derivative of the $n$-th coefficient in both cases as long as $g_0(x)=1$. It is interesting however, if there are any useful applications of such a general formula in other cases besides the binomial one.


\begin{thebibliography}{100}
\bibitem[Dlt]{Dlt}
G.-C. Rota, D. Kahaner, and A. Odlyzko, "Finite Operator Calculus," Journal of Mathematical Analysis and its Applications, vol. 42, no. 3, June 1973. Reprinted in the book with the same title, Academic Press, New York, 1975.

\bibitem[Gdsl]{Gdsl}	
C. D. Godsil, Algebraic Combinatorics, Chapman \& Hall, New York, 1993. Chapter 9.

\bibitem[LB]{LB}
D. Krotkov, ''Taking the logarithm of binomial type sequences: linear approach''\\
https://arxiv.org/abs/1910.07100

\end{thebibliography}
\end{document}